\documentclass[a4paper]{amsart}
\usepackage{amsmath}
\usepackage{amsfonts}
\usepackage{amsthm}
\usepackage{amssymb}
\usepackage[english]{babel}
\usepackage[latin1]{inputenc}
\usepackage[labelsep=space]{caption}
\usepackage{graphicx,color}
\usepackage{enumerate}

\newcommand{\inte}{{\mathop{\mathrm{int}\,}}}
\newcommand{\cl}{{\mathop{\mathrm{cl}\,}}}

\newcommand{\conv}{{\mathop{\mathrm{conv}\,}}}

\newcommand{\supp}{{\mathop{\mathrm{supp}\,}}}

\newcommand{\pa}{{\partial}}

\newcommand{\Real}{\mathbb{R}}

\newcommand{\ee}{{\epsilon}}

\newcommand{\de}{{\delta}}

\newcommand{\la}{{\lambda}}

\newcommand{\rotpih}{\phi_{\pi/2}}

\newtheorem{theorem}{Theorem}[section]
\newtheorem{corollary}[theorem]{Corollary}
\newtheorem{lemma}[theorem]{Lemma}
\newtheorem{proposition}[theorem]{Proposition}

\theoremstyle{definition}

\numberwithin{equation}{section}

\begin{document}

\date{\today}
\title{Covariogram of non-convex sets}

\author{Carlo Benassi, Gabriele Bianchi and Giuliana D'Ercole }
\address{Dipartimento di Matematica,
Universit\`a di Modena e Reggio Emilia, via Campi 213/B,
Modena, Italy  I-41100} \email{benassi.carlo@unimore.it}\email{giulianadercole@libero.it}
\address{Dipartimento di Matematica, Universit\`a di Firenze, Viale Morgagni 67/A,
Firenze, Italy I-50134} \email{gabriele.bianchi@unifi.it}

\begin{abstract}
\noindent The covariogram of a compact set $A\subset\mathbb{R}^n$
is the function that to each $x\in \mathbb{R}^n$ associates
the volume of $A\cap (A+x)$. Recently it has been proved that the covariogram determines any planar
convex body, in the class of all convex bodies. We extend the class of sets in which a planar convex body
is determined by its covariogram. 
Moreover, we prove that there is no pair of  non-congruent planar polyominoes consisting of less than $9$ points that have equal discrete covariogram.
\end{abstract}
\maketitle

\section{Introduction}
Let $A$ be a compact set in $\Real^n$, $n\geq2$, and let $\la_n$ stand for the $n$-dimensional Lebesgue measure. The \emph{covariogram} $g_A$ of  $A$ is the function on $\Real^n$ defined by
\begin{equation}
g_A(x) := \la_n(A \cap (A+x)), \qquad x \in \Real^n.
\end{equation}
This function is also called \emph{set covariance}, and it  coincides with the  \emph{autocorrelation} of the characteristic function $\mathbf{1}_A$ of $A$, i.e. $g_A =\mathbf{1}_A\ast \mathbf{1}_{-A}$.
The covariogram $g_A$ is clearly unchanged by a translation or a reflection 
of $A$. (The term \emph{reflection} will always mean reflection 
in a point.) A \emph{convex body} in $\Real^n$ is a convex compact set with non-empty interior. In 1986 Matheron~\cite[p.~20]{Ma2}  asked the following question and conjectured a positive answer for the case $n=2$. (The same question was also asked independently by Adler and Pyke \cite{AP} in the  probabilistic terms expressed by Problem~P\ref{distribution_XY} below.)
\smallskip

\textbf{Covariogram problem.} \emph{Does $g_K$ determine a
convex body $K$, among all convex bodies, up to translations and
reflections?}
\smallskip

The conjecture for $n=2$ has been completely settled only very recently, by  Averkov and Bianchi~\cite{AB+07}.
\begin{theorem}[\cite{AB+07}] \label{MatConjConf}
Every planar convex body $K$ is determined within all planar convex bodies by its covariogram, up to translations and reflections.
\end{theorem}
 
See \cite{AB+07} for further information on the covariogram problem. In general, the convexity of $K$  is needed in Theorem~\ref{MatConjConf}, since there exist pairs  of  non-convex non-congruent (i.e. there exists no isometry mapping one into the other) planar polygons with equal covariogram; see Gardner, Gronchi and Zong~\cite{GGZ05} and Fig.~\ref{figronchi}.
We prove some results which  extend the class of bodies in which a convex body is determined by its covariogram. The main results of this type are the following two. 
Let $\mathcal{A}$ denote the class of planar regular (i.e. equal to the closure of their interior) compact sets whose interior has at most two  components.

\begin{theorem}\label{nonconvessi3}
If $A\in \mathcal{A}$ and $g_A=g_K$, for some convex body $K$
in $\mathbb{R}^2$, then $A$ is convex.
\end{theorem}

Let $\mathcal{B}$ denote the class of planar compact sets whose boundary consists of a finite number of closed disjoint simple polygonal curves (each one with finitely many edges). The class $\mathcal{B}$ contains each set which is finite union of disjoint polygons, as well as sets that can be written as $P\setminus Q$, with $P$ and $Q$ polygons and $Q\subset\inte P$.

\begin{theorem}\label{poli}
If $A\in \mathcal{B}$ and $g_A=g_K$, for some convex body $K$
in $\mathbb{R}^2$, then $A$ is convex.
\end{theorem}

The  previous theorems clearly imply that a planar convex body is determined by its covariogram  both in the class $\mathcal{A}$ and $\mathcal{B}$. 
It is known that the covariogram problem is equivalent to any of the following problems (see~\cite{AB+07} for a detailed explanation of each problem with references to the relevant literature):
\begin{enumerate}[{P}1]
\item\label{chord_lengths} Determine a convex  body $K$ by the knowledge, for each unit vector $u$ in $\Real^n$, of the distribution of the lengths of the chords of $K$ parallel to $u$.
\item\label{distribution_XY} 
Determine a convex  body $K$ by the  distribution of $X-Y$, where $X$ and $Y$ are independent random variables uniformly distributed over $K$.
\item\label{phase_retrieval} Determine  the characteristic function $\mathbf{1}_K$ of a convex body $K$ from  the modulus of its Fourier transform $\widehat{\mathbf{1}_K}$.
\end{enumerate}
Thus the previous theorems imply a positive answer to Problems~P\ref{chord_lengths} and P\ref{distribution_XY} both in the class $\mathcal{A}$ and $\mathcal{B}$, and a positive answer to Problem~P\ref{phase_retrieval} in the class of characteristic functions of sets in $\mathcal{A}$ or in $\mathcal{B}$.

Propositions~\ref{insdiff}, \ref{nonconvessi1}, \ref{nonconvessi2}, \ref{cond_volume} and Corollary~\ref{nonconvessi4}, all contained in Section~\ref{conv-non-conv}, are other results in the spirit of Theorems~\ref{nonconvessi3} and \ref{poli}. 

Some aspects of the covariogram problem are of combinatorial
nature. Two finite subsets $A$ and $B$ of $\mathbb{R}^n$ are said to be \emph{homometric} if $A\cap (A+x)$ and $B\cap (B+x)$ have equal cardinality for each $x\in \mathbb{R}^n$, or, equivalently, if the sets of vector differences $\{x-y:x,y\in A\}$ and $\{x-y:x,y\in B\}$ are identical counting multiplicities.
One problem consists in determining all the sets which are homometric to a given set.
We refer to \cite{RS}, \cite{LRH} and \cite{DGN02} for a complete algebraic solution of this problem for subsets of the real line.

A \emph{polyomino} is a finite subset $A$ of $\mathbb{Z}^n$ such that the union $A+[0,1]^n$ of lattice unit cubes has connected interior. A polyomino $A$ is \emph{convex} if $A=(\conv A)\cap \mathbb{Z}^n$. We shall refer to the set $A+[0,1]^n$ (itself called a polyomino by many authors) as the \emph{animal} of the polyomino $A$. The non-convex polygons with equal covariogram presented in \cite{GGZ05} are the animals of two homometric convex polyominoes made of $15$ points. We are interested in finding a similar example with minimal number of unit squares. Since  two polyominoes are homometric if and only if the associated animals have the same covariogram (see Lemma~\ref{poly_animals}) we are interested in finding pairs of homometric polyominoes with minimal cardinality. 
We exhibit a pair of non-congruent homometric polyominoes made of $9$ points, and 
we prove that this example is minimal. 
\begin{theorem}\label{polyominoes} The minimum number $d$ such that
there exists a pair of non-congruent homometric polyominoes in $\mathbb{Z}^2$
made of $d$ points is $9$.
\end{theorem}
In terms of animals, this theorem proves that two non-congruent animals made of less than nine unit squares cannot have the same covariogram.

%
%

\section{Definitions, notations and preliminaries}\label{notations}

As usual, $\mathbb{S}^{n-1}$ denotes the unit sphere in Euclidean $n$-space
$\mathbb{R}^n$. For $x, y \in \mathbb{R}^n$, $\|x\|$ denotes the Euclidean norm of $x$, $x\cdot y$ denotes scalar product, while $[x,y]$ denotes the closed line segment with endpoints $x$ and $y$. For $\ee>0$ the symbol $B(x,\ee)$ denotes the open ball centred at $x$ and with radius $\ee$.
If $u\in\mathbb{S}^{n-1}$, we denote by $u^{\perp}$ the
$(n-1)$-dimensional subspace orthogonal to $u$ and by $l_u$ the line parallel
to $u$ containing the origin $o$. The symbol $\rotpih$ denotes counterclockwise rotation by $\pi/2$ about the origin in $\Real^2$.
We write $\lambda_k$ for $k$-dimensional Lebesgue measure in $\mathbb{R}^n$, where $k=1,\ldots,n$,
and where we identify $\lambda_k$ with $k$-dimensional Hausdorff measure. 

If $A$ is a set, we denote by  $|A|$, $\cl A$, $\inte A$ and $\conv A$ the \emph{cardinality}, \emph{closure},
\emph{interior}, and \emph{convex hull} of $A$, respectively. The notation for
the usual orthogonal \emph{projection} of $A$ on a subspace $S$ is $A|S$.
A closed set $A$ is said to be \emph{regular} if it coincides with the closure of its interior.
If $A$ and $B$ are subsets of $\mathbb{R}^n$, their \emph{Minkowski sum} is
\[
A+B=\{a+b: a\in A, b\in B\}.
\]
In particular, if $x\in\mathbb{R}^n$, then $A+x$ denotes the translate of $A$ by
$x$. We also write $DA=A+(-A)$ for the \emph{difference set} of $A$. 

The \emph{support function} of a compact set $A\subset\Real^n$ is defined, for $x\in\Real^n$, by 
\[
h_A(x)=\sup \{x\cdot y : y\in A\},
\]
while the \emph{width of $A$ in direction $u\in \mathbb{S}^{n-1}$} is defined by
\[
w(A,u)=h_A(u)+h_A(-u).
\]
The linearity of the support function with respect to Minkowski addition implies $w(A,u)=(1/2)w(DA,u)$. 

It is well known that, when $A$ is compact, $g_A$ is continuous and 
\[
 \supp g_A=DA, 
\]
where $\supp f$ denotes the \emph{support} of the function $f$.

Given $u\in \mathbb{S}^{n-1}$ and a compact set $A\subset\Real^n$, the \emph{exposed face of $A$ in direction $u$} is $F(A,u)=\{x\in A : x\cdot u=h_A(u)\}$. \cite[Th.~1.7.5(c)]{Schn} proves that, for a convex body $A$  and $u\in \mathbb{S}^{n-1}$, 
\begin{equation}\label{facce_corpodiff}
F(DA,u)=F(A,u)-F(A,-u),
\end{equation}
and is it not difficult to see that the previous formula is valid also  for any compact set $A$.




\section{Comparison between covariograms of convex and non-convex
sets}\label{conv-non-conv}

Proposition~\ref{insdiff} exploits the convexity of $\supp g_K=D\,K$ when $K$ is  a convex body.

\begin{lemma}\label{inviluppoconvesso}
Given a compact set $A\subset\mathbb{R}^n$, one has $\conv(DA)=D(\conv A)$.
\end{lemma}

\begin{proof} It suffices to prove $h_{\conv(DA)}=h_{D(\conv A)}$.  This identity is a consequence of the linearity of the support function with respect to Minkowski addition, and of the identity $h_B=h_{\conv B}$, valid for each compact set $B$. Indeed, one has
\[
h_{\conv (DA)}=h_{DA}=h_{A-A}=h_A+h_{-A}=h_{\conv A}+h_{-\conv A}=h_{D(\conv A)}.
\]
\end{proof}

\begin{proposition}\label{insdiff}
Let $A$ be a regular compact set of $\mathbb{R}^n$. If either $\supp g_A$ is not convex, or $D\,A\neq D(\conv A)$, then $g_A\neq g_K$, for each convex body $K\subset\mathbb{R}^n$.
\end{proposition}

\begin{proof} If $A$ has the same covariogram as a convex
body, then $DA$ is convex, and so $DA=\conv(DA)$. By Lemma
\ref{inviluppoconvesso}, this implies $DA=D(\conv A)$. 
\end{proof}

It is well known (see \cite[p.86]{Ma1}) that, when $K$ is a convex body and $u\in \mathbb{S}^1$,
\begin{equation}\label{derivata_convessi2}
-\frac{\partial g_K}{\partial u}(0)=\la_{n-1}(K\mathbin|u^\perp).
\end{equation}
When $n=2$, since $\la_{1}(K\mathbin|u^\perp)=w(K,\rotpih u)=({1}/{2})w(\supp g_K,\rotpih u)$, the formula becomes
\begin{equation}\label{derivata_convessi}
-\frac{\partial g_K}{\partial u}(0)=\frac{1}{2}w(\supp g_K,\rotpih u).
\end{equation}
Propositions~\ref{nonconvessi1} and~\ref{nonconvessi2} and Theorem~\ref{nonconvessi3} exploit \eqref{derivata_convessi} to give conditions on a planar regular compact set $A$ that imply $g_A\neq g_K$ for every convex body $K$. We begin with two lemmas.

\begin{lemma}\label{mono}
Let $D$ be a bounded Lebesgue-measurable subset of $\mathbb{R}$
such that $\lambda_1(D)>0$. Then 
$$
\liminf_{h\rightarrow
0}\frac{\lambda_1(D\setminus(D+h))}{h}\geqslant 1.
$$
\end{lemma}

\begin{proof} As $\lambda_1(D)>0$, there exists a point
$x_0$ of density for $D$, i.e.\  a point such that for every
positive $\epsilon$ there exists $\bar{h}(\epsilon)>0$ such that
for every  $h\in(0, \bar{h}(\ee))$ we have
$$
\lambda_1\left(\left(x_0-h/2,x_0+h/2\right)\cap
D\right)\geqslant(1-\epsilon)h
$$ 
(see \cite [Cor. 6.26]{Coh}).
Choose $\ee>0$ and let $h\in(0, \bar{h}(\ee))$.
Let $D_0=\left(x_0-h/2,x_0+h/2\right)\cap
D$ and consider the sequence $(B_k)$, where
$$
B_k=\left(D_0\cap (D+h)\cap ... \cap (D+(k-1)h)\right)\setminus (D+kh).
$$
We have $D_0=\bigcup_{i=1}^\infty B_i$. In fact, the inclusion $D_0\supseteq\bigcup_{i=1}^\infty B_i$ is trivial; while if
$x\in D_0$, then there exists $k$ such that $x\notin D+kh$, since $D$ is bounded. Let $\bar{k}$ be the smallest integer such that $x\notin D+\bar{k}h$. Then $x\in B_{\bar{k}}\subset\bigcup_{i=1}^\infty B_i$.
If $i\neq j$, then $B_i\cap B_j=\emptyset$.  In fact, when $i<j$ we have $B_i\cap(D+ih)=\emptyset$ and $B_j\subset(D+ih)$. 
It follows that 
\begin{equation}\label{lemma1_1}
\sum_{i=1}^\infty\lambda_1(B_i)=\lambda_1(D_0)\geqslant (1-\epsilon)h.
\end{equation}
If $i\neq j$, then  $(B_i-(i-1)h)\cap (B_j-(j-1)h)=\emptyset$. In fact, $B_i\subset\left(x_0-h/2,x_0+h/2\right)$ implies  $B_i-(i-1)h\subset I_i:=\left(x_0-(i-1)h-h/2, x_0-(i-1)h+h/2\right)$ and $I_i\cap I_j=\emptyset$.
Let us also remark that $B_i\subset(D+(i-1)h)\setminus(D+ih)$ and so $B_i-(i-1)h\subset D\setminus(D+h)$.
Thus, 
\begin{equation}\label{lemma1_2}
\sum_{i=1}^\infty\lambda_1(B_i)=\sum_{i=1}^\infty\lambda_1(B_i-(i-1)h)\leqslant \lambda_1(D\setminus(D+h)).
\end{equation}

The inequalities \eqref{lemma1_1} and \eqref{lemma1_2} imply the statement. \end{proof}

\begin{lemma}\label{mono2}
Let $D$ be a bounded Lebesgue-measurable subset of $\mathbb{R}$ such that there exist
$2r$ points $a_1, b_1,\dots, a_r, b_r
\in\mathbb{R}\setminus D$, with $a_1<b_1<a_2<b_2<...<a_r<b_r$, for
which $\lambda_1(D \cap [a_i, b_i])>0$, $i=1,\ldots, r$, and
$\lambda_1(D \cap [b_i, a_{i+1}])=0$, $i=1,\ldots, r-1$. Then
$$\liminf_{h\rightarrow 0}\frac{\lambda_1(D\setminus(D+h))}{h}\geqslant r.$$
\end{lemma}

\begin{proof} Let $D_i=D\cap[a_i, b_i]$ and let $h$ satisfy $|h|<\min_{i\neq j}|b_i-a_{j}|$. Observe that if $i\neq j$, then $D_i\cap(D_j+h)=\emptyset$. Therefore
$$
g_D(h)=\lambda_1(D\cap(D+h))=\lambda_1\left(\bigcup_{i=1}^r(D_i\cap(D_i+h)\right)=\sum_{i=1}^r\lambda_1((D_i\cap(D_i+h)).
$$
Thus, 
$$
\lambda_1(D\setminus(D+h))=\lambda_1(D)-g_D(h)=\lambda_1(D)-\sum_{i=1}^r\lambda_1((D_i\cap(D_i+h))=
$$
$$
=\lambda_1(D)-\sum_{i=1}^r(\lambda_1(D_i)-\lambda_1(D_i\setminus(D_i+h)))=\sum_{i=1}^r\lambda_1(D_i\setminus(D_i+h)).
$$
The statement follows by applying Lemma \ref{mono} to each set $D_i$. 
\end{proof}

Let $A$ be a regular compact set of $\mathbb{R}^2$, $u\in\mathbb{S}^1$ and $y\in u^\perp$. Let us set
\[
f_{A,u}(y)=\liminf_{h\rightarrow 0}\frac{\lambda_1((A\setminus(A+hu))\cap(y+l_u))}{h}.
\]

\begin{proposition}\label{nonconvessi1}
Let $A\subset\mathbb{R}^2$ be a regular compact set for which there exists a direction $u\in\mathbb{S}^1$ such that
\begin{enumerate}
\item\label{prop1_i} $\lambda_1(A\cap(y+l_u))>0$ for $\la_1$-a.e. $y\in \conv(A|u^{\perp})$ and
\item\label{prop1_ii} $\lambda_1(\{y\in u^\perp:f_{A,u}(y)\geqslant2\})>0$.
\end{enumerate}
\noindent Then $g_A\neq g_K$, for every convex body $K$ in $\mathbb{R}^2$.
\end{proposition}

\begin{proof} 
If $(\partial g_A/\partial
u)(0)$ does not exist, the statement follows by \eqref{derivata_convessi}. Otherwise
\begin{equation}\label{misura_fubini}\begin{aligned}
-\frac{\partial g_A}{\partial u}(0)=&
\lim_{h\rightarrow
0}\frac{\lambda_2(A\setminus(A+hu))}{h}\\
=&\lim_{h\rightarrow
0}\int_{\conv(A|u^\bot)}\frac{\lambda_1((A\setminus(A+hu))\cap(y+l_u))}{h}
d\la_1(y),
\end{aligned}\end{equation}
since $A|u^\bot=\conv(A|u^\bot)$, by Assumption~\eqref{prop1_i} and the fact that $A\mathbin|u^\perp$ is closed.
Let $A^{(2)}=\{y\in u^\perp:f_{A,u}(y)\geqslant2\}$. By virtue of Fatou's lemma we have
\begin{multline*}
-\frac{\partial g_A}{\partial u}(0)\geqslant
\int_{\conv(A|u^\bot)\setminus A^{(2)}}
\liminf_{h\rightarrow0}\frac{\lambda_1((A\setminus(A+hu))\cap(y+l_u))}{h}d\la_1(y)+
\\
+\int_{A^{(2)}}
\liminf_{h\rightarrow0}\frac{\lambda_1((A\setminus(A+hu))\cap(y+l_u))}{h}d\la_1(y).
\end{multline*}
Thus, by Lemmas~\ref{mono} and~\ref{mono2} and Assumptions~\eqref{prop1_i} and~\eqref{prop1_ii}, we have 
\[
-\frac{\partial g_A}{\partial u}(0)
\geqslant 
\int_{\conv(A|u^\bot)\setminus A^{(2)}} d\la_1(y)+\int_{A^{(2)}} 2\ d\la_1(y)
>\lambda_1(\conv(A|u^{\bot})).
\]
Since $\lambda_1(\conv(A|u^{\bot}))=w(A,\rotpih u)$, and $w(A,\rotpih u)=(1/2)w(\supp g_A,\rotpih u)$ (because $\supp g_A=D\,A$), we have
\[
-\frac{\partial g_A}{\partial u}(0)
>\frac{1}{2}w(\supp g_A,\rotpih u).
\]
This inequality and \eqref{derivata_convessi} imply $g_A\neq g_K$,
for every convex body $K$ in $\mathbb{R}^2$. 
\end{proof}

\begin{corollary}\label{nonconvessi4}
Let $A\subset\mathbb{R}^2$ be a regular compact set such that $\inte A$ has finitely many  components. Assume that there exist $u\in\mathbb{S}^1$ and $a_1$, $a_2\in\inte A$ such that 
$A\mathbin|u^\perp$ is a segment and $[a_1,a_2]$ is parallel to $u$ and meets $\Real^2\setminus A$. Then $g_A\neq g_K$, for every convex body $K$ in $\mathbb{R}^2$.
\end{corollary}

\begin{proof} It suffices to prove that  the assumptions of  Proposition~\ref{nonconvessi1} are satisfied.
The assumptions of the corollary imply that $(\inte A)\mathbin| u^\perp$ consists of finitely many intervals and $A\mathbin| u^\perp$ is a regular closed set. Thus $\la_1(A\mathbin|u^\perp\setminus(\inte A)\mathbin| u^\perp)=0$. Since,  $\lambda_1(A\cap(y+l_u))$ is positive when $y\in (\inte A)\mathbin| u^\perp$, Assumption~\eqref{prop1_i} of Proposition~\ref{nonconvessi1} is satisfied.

Let $b\in[a_1,a_2]\setminus A$, and let $\ee>0$ be such that $B(a_i,\ee)\subset A$, $i=1,2$, and $B(b,\ee)\subset \Real^2\setminus A$. If $y\in u^\perp\cap B(a_1|u^{\perp},\ee)$, then  $A\cap(y+l_u)$
contains two closed non-degenerate intervals separated by a non-degenerate interval contained in $\Real^2\setminus A$. Thus, by Lemma~\ref{mono2}, $u^\perp\cap B(a_1|u^{\perp},\ee)\subset\{y\in u^\perp:f_{A,u}(y)\geqslant2\}$ and Assumption~\eqref{prop1_ii} of Proposition~\ref{nonconvessi1} is satisfied. 
\end{proof}

%

\begin{proof}[Proof of Theorem~\ref{nonconvessi3}] We argue by contradiction. Assume that
$\inte A$ has two  components $A_1$ and $A_2$, and let 
$a_1$ and $a_2$  belong respectively to  $A_1$ and $A_2$, and be such that $[a_1,a_2]$ meets $\Real^2\setminus A$. 
Let $u$ be the direction of the segment $[a_1,a_2]$. The set $(\inte A)\mathbin|u^\perp$ is an interval, because $A_1\mathbin|u^\perp$ and $A_2\mathbin|u^\perp$ are intervals and, by the definition of $u$, their intersection is non-empty. Thus, by Corollary~\ref{nonconvessi4}, $g_A\neq g_K$, for every convex body $K$ in $\mathbb{R}^2$. 

When $\inte A$ has only one  component the proof is
similar. Since $A$ is supposed to be non-convex, one can choose $a_1, a_2\in
\inte A$ so that $[a_1, a_2]\cap (\mathbb{R}^2\setminus A)\neq\emptyset$. Again Corollary~\ref{nonconvessi4} gives a contradiction. 
\end{proof}

When, for each line $l$ parallel to $u\in \mathbb{S}^1$, the section $A\cap l$ is, up to a set of measure zero, the union of closed segments, Proposition~\ref{nonconvessi1} can be made more precise. 

\begin{proposition}\label{nonconvessi2}
Let $A\subset\mathbb{R}^2$ be a regular compact set. Assume that there
exists $u\in \mathbb{S}^1$ such that, for each $y\in u^\perp$, the section $A\cap(y+l_u)$ consists, up to a set of $\la_1$-measure zero, of a finite or infinite  number $N(y)$  of closed disjoint segments. 
If
\begin{equation}\label{hp}
\sum_{i=0}^\infty(i-1)
\lambda_1\{y\in \conv(A\mathbin|u^{\bot}):N(y)=i\}\neq0,
\end{equation}
then $g_A\neq g_K$, for
every convex body $K$ in $\mathbb{R}^2$.
\end{proposition}

\begin{proof} Let us first prove that, for each $h\neq0$,
\begin{equation}\label{eqn}
\frac{\lambda_1((A\setminus (A+hu))\cap(y+l_u))}{h} \leqslant
N(y).
\end{equation} 
For brevity, let $l=y+l_u$. If $A\cap l=\bigcup_{i=1}^{N(y)}[a_i, b_i]$, up to a set of measure zero, then
\begin{align*}
\lambda_1((A\setminus (A+hu))\cap l)&=
\lambda_1(A\cap l)-\lambda_1(A\cap (A+hu)\cap l)\\
&=\sum_{i=1}^{N(y)}\lambda_1([a_i, b_i])-
\sum_{i=1}^{N(y)}\lambda_1([a_i, b_i]\cap[a_i+h, b_i+h])\\
&\quad\quad
-\sum_{i\neq j,\  i,j=1}
^{N(y)}\lambda_1([a_i, b_i]\cap[a_j+h, b_j+h])\\
&\leqslant
\sum_{i=1}^{N(y)}\lambda_1\left([a_i, b_i]\setminus [a_i+h,
b_i+h]\right).
\end{align*}
Since each summand in the last sum of the previous formula is less than or equal to $h$, we have \eqref{eqn}. 

Lemma~\ref{mono2} implies that, for every integer  $r$ with $r\leqslant N(y)$, we have 
\begin{equation}\label{eqn2}
\liminf_{h\to0}\frac{\lambda_1((A\setminus (A+hu))\cap(y+l_u))}{h} 
\geqslant r. 
\end{equation}
Formulas~\eqref{eqn}  and~\eqref{eqn2} imply 
\begin{equation*}\label{eqn3}
\lim_{h\to0}\frac{\lambda_1((A\setminus (A+hu))\cap(y+l_u))}{h} =N(y).
\end{equation*}
If $\int_{\conv(A\mathbin|u^\perp)} N(y)\ d\la_1(y)=+\infty$, then 
$-(\partial g_A/\partial u)(0)=+\infty$, by \eqref{misura_fubini} and Fatou's Lemma. In this case \eqref{derivata_convessi} implies $g_A\neq g_K$ for every convex body $K$. If the previous integral is finite, 
then we may apply the Lebesgue dominated convergence Theorem to the last integral in \eqref{misura_fubini}, and we have
\[
-\frac{\partial g_A}{\partial u}(0)=\int_{\conv(A|u^\bot)}N(y)\ d\la_1(y).
\]
Since
\begin{align*}
\int_{\conv(A|u^\bot)}N(y)\ d\la_1(y)&=\sum_{i=0}^\infty\int_{\{y\in \conv(A|u^\bot):N(y)=i\}}i\ d\la_1(y)\\
&=\sum_{i=0}^\infty i\  \lambda_1\{y\in \conv(A|u^\bot):N(y)=i\},
\end{align*}
and $({1}/{2})w(\supp g_A,\rotpih u)=w(A, \rotpih u)=\sum_{i=0}^\infty \lambda_1\{y\in \conv(A|u^\bot):N(y)=i\}$, \eqref{hp} implies
\[
-\frac{\partial g_A}{\partial u}(0)\neq \frac{1}{2}w(\supp g_A,\rotpih u).
\]
Again, \eqref{derivata_convessi} implies $g_A\neq g_K$ for every convex body $K$.
\end{proof}

The next result is valid for sets of any dimension. The covariogram $g_A$ provides both $\la_n(A)=g_A(0)$ and $\la_n(DA)=\la_n(\supp g_A)$. Since when $A$ is convex $\la_n(A)$ and $\la_n(D\,A)$ are related by the Rogers-Shephard and the Brunn-Minkowski  inequalities, we obtain some conditions on $g_A$ which are necessary for $A$ to be convex.

\begin{proposition}\label{cond_volume}
Let $A\subset\mathbb{R}^n$ be a regular compact set. If $A$ is convex, then 
\begin{gather}
{\binom{2n}{n}}^{-1}\lambda_n(\supp g_A)\leqslant
g_A(0)\leqslant 2^{-n}\lambda_n(\supp g_A)\label{necess1}
\intertext{and, for each $u\in \mathbb{S}^{n-1}$,}
{\binom{2n-2}{n-1}}^{-1}
\lambda_{n-1}(\supp g_A\mathbin|u^\bot)\leqslant
-\frac{\partial g_A} {\partial
u}(0)\leqslant2^{1-n}\lambda_{n-1}(\supp g_A\mathbin|u^\bot).\label{necess2}
\end{gather}
\end{proposition}

\begin{proof} 
%
The Rogers-Shepard and the Brunn-Minkowski inequalities (see~\cite[Th.~7.3.1]{Schn}) state that, when $A\subset\Real^n$ is convex, we have
\[
{\binom{2n}{n}}^{-1}\lambda_n(DA)\leqslant
\la_n(A)\leqslant 2^{-n}\lambda_n(DA).
\]
Thus \eqref{necess1} is an immediate consequence of the previous inequalities and of the identities $g_A(0)=\la_n(A)$ and $DA=\supp g_A$. The same inequalities, applied to the $(n-1)$-dimensional convex body $A\mathbin|u^\perp$, give
\[
{\binom{2n-2}{n-1}}^{-1}
\lambda_{n-1}(D(A\mathbin|u^\perp))\leqslant
\la_{n-1}(A\mathbin|u^\perp)\leqslant2^{1-n}\lambda_{n-1}(D(A\mathbin|u^\perp)).
\]
The identity $D(A\mathbin|u^\perp)=(DA)\mathbin|u^\perp=\supp g_A \mathbin|u^\perp$ and \eqref{derivata_convessi2} imply \eqref{necess2}.
%
\end{proof}

%
%
%

\begin{figure}[h]
\centering
\includegraphics[scale=0.52]{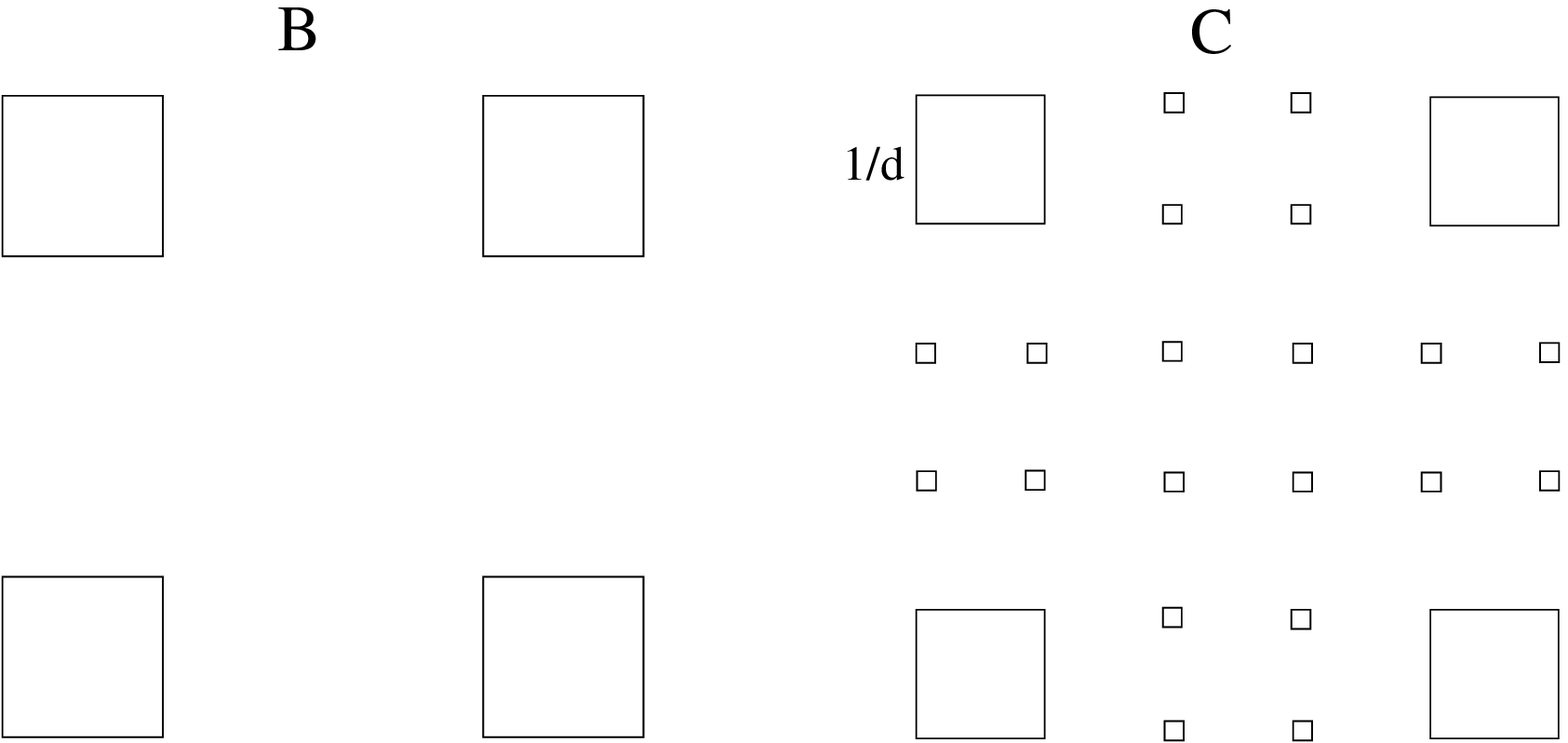}
\caption{}\label{figura_Ab}
\end{figure}

\begin{figure}[h]
\centering
\includegraphics[scale=0.6]{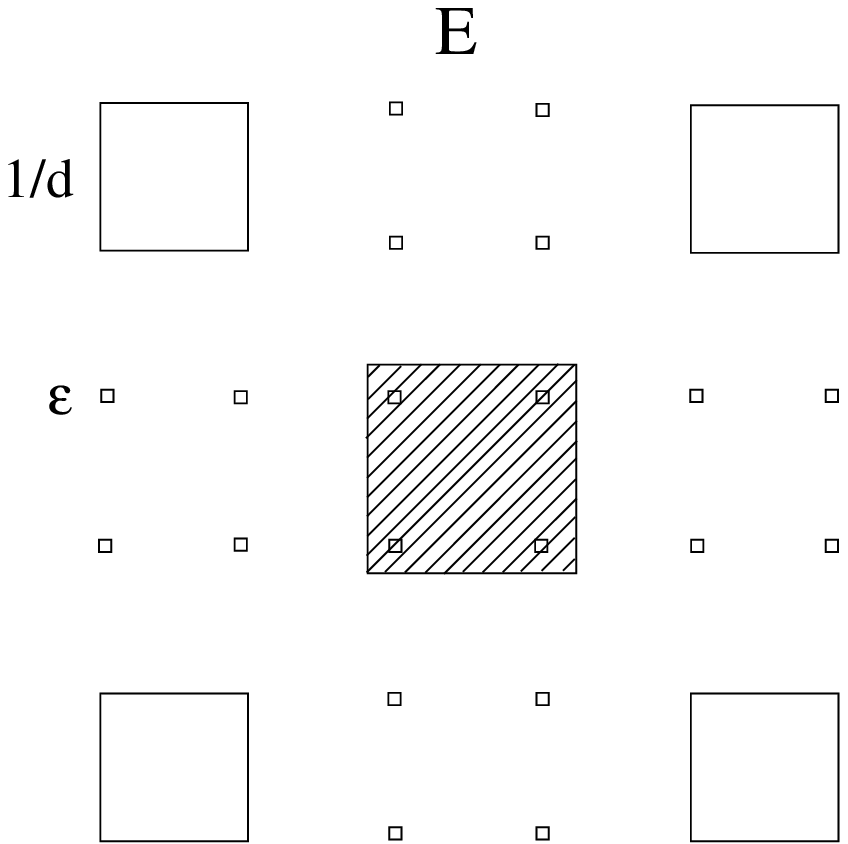}
\caption{}\label{figura_B}
\end{figure}

In order to critically discuss the  previous results, let us present some examples (see Figures~\ref{figura_Ab} and \ref{figura_B}). 
%
Let $Q=[0,1]^2$. The set $B$ is obtained by placing four squares of edge $1/4$ inside and in the corner of  $Q$, so that $\conv B=Q$. To prove that $g_B$ differs from the covariogram of any convex body one cannot use  Proposition~\ref{nonconvessi2}, because \eqref{hp} is false, but one can use Proposition~\ref{insdiff}, since $D\,B\neq D\,Q$ is not convex. 
%
The set $C$ is constructed as follows. Divide $Q$ in $d^2$ equal squares.
We obtain a grid of $(d+1)^2$ points. The body $C$ is the subset
of $Q$ which is the union of the four squares of edge $1/d$
touching the four vertex of $Q$ and of $(d+1)^2-16$ little
squares of edge 
$\ee=(1-4/d)/((d+1)^2-16)$ contained in $Q$ and
containing the points of the grid outside the four
squares already considered. It results that $C$ does not satisfy
condition \eqref{hp} in Proposition \ref{nonconvessi2} and, moreover, 
$D\,C=D\,Q$ is convex. In this case, what proves that   $g_C$ differs from the covariogram of a convex body when $d$ is large is Proposition~\ref{cond_volume}, since \eqref{necess1} is not satisfied by $C$ (because $1/6 \la_1(D\,Q)>g_C(0)$ when $d$ is large). 
Choose $\epsilon$ so that $0<\epsilon<(1-4/d)/((d+1)^2-16)$. The set $E$ (see Fig.
\ref{figura_B}) is constructed by adding
another square in the centre of $C$ of edge
$1-{4}/{d}-((d+1)^2-16)\epsilon$ (actually, a little bit
longer than this, to compensate for the little squares  included in this central square which disappear so that \eqref{hp} does not hold).  The set $E$ does not satisfy  (\ref{hp}), we have $DE$ convex and, when $d$ is large and $\epsilon$ is small even \eqref{necess1} and \eqref{necess2} are satisfied.
%
%
The fact that a set like $E$ does not have the covariogram equal to that of a convex body is a consequence of Theorem~\ref{poli}.

%

In order to prove Theorem~\ref{poli} we need a lemma computing a second order distributional derivative of $g_A$. These computations are made in  \cite[Lemma~4.2]{Bi06}  when $A$ is a convex polytope in $\Real^n$, and can be repeated, almost without any change, also when $A\in\mathcal{B}$. 
Let $C^\infty_0(\mathbb{R}^2)$ denote the class of infinitely differentiable functions on $\mathbb{R}^2$ with compact support. We recall that $|A|$ denotes the cardinality of $A$.

\begin{lemma}\label{utile}
Let $A\in\mathcal{B}$ and $F_1,\ldots,F_m$ be the edges of the polygons which constitute $A$. Let $\nu_i$, $i=1,\ldots, m$, be the unit outer
normal vector of $A$ at $F_i$, $w\in \mathbb{S}^1$,  $I_p=\{(i,j):\text{$F_i$ is
parallel to  $F_j$}\}$ and $I_{np}=\{(i,j):\text{$F_i$ is not parallel to  $F_j$}\}$.
Then, for $f\in C^\infty_0(\mathbb{R}^2)$, we have
\begin{multline}\label{formulone}
-\frac{\partial^2 g_A}{\partial w^2}\left(f\right)=
\sum_{(i,j)\in I_{np}}
\frac{w\cdot\nu_i\ w\cdot\nu_j}{\sqrt{1-(\nu_i\cdot\nu_j)^2}}
	\int_{\mathbb{R}^2} |F_i\cap(F_j+z)|\,f(z)\,d\la_2(z) +\\
+\sum_{(i,j)\in I_{p}} 
w\cdot\nu_i\ w\cdot\nu_j
\int_{F_i-F_j}\lambda_{1}(F_i\cap(F_j+z))\,f(z)\,d\la_1(z).
\end{multline}

Both sums in the right hand side of (\ref{formulone}) are uniquely determined by $g_A$.
\end{lemma}
\begin{proof}
The definition of derivative in the sense of distributions implies $(\pa 1_A / \pa w) (f)=-\int_{A}\pa f(x)/\pa w\, dx$. Thus, by the Divergence Theorem, we have
\begin{equation*}
\frac{\pa 1_A}{\pa w} (f)=-\sum_{i=1}^m w\cdot\nu_i\de_{F_i} (f),
\end{equation*}
where $\de_{F_i} (f)=\int_{F_i}f(x)d\la_{1}(x)$.
Since $g_A=1_A\ast1_{-A}$,  we can write
\begin{equation*}\label{somma_dirac}
\frac{\pa^2 g_A}{\pa w^2} (f)=
\left(\frac{\pa1_A}{\pa w}\ast\frac{\pa1_{-A}}{\pa w}\right) (f)
=-\sum_{i,j=1}^m w\cdot\nu_i\ w\cdot\nu_j\  (\de_{F_i}\ast\de_{-F_j}) (f).
\end{equation*}
A direct computation (see \cite[Lemma 4.2]{Bi06} for the details) proves 
\begin{gather*}
(\de_{F_i}\ast\de_{-F_j}) (f)=
\int_{F_i-F_j}\lambda_{1}(F_i\cap(F_j+z))\,f(z)\,d\la_1(z)
\intertext{when $F_i$ and $F_j$ are parallel, and}
(\de_{F_i}\ast\de_{-F_j}) (f)=
(1-{(\nu_i\cdot \nu_j)}^2)^{-1/2}
\int_{\mathbb{R}^2} |F_i\cap(F_j+z)|\,f(z)\,d\la_2(z)
\end{gather*}
when $F_i$ and $F_j$ are not parallel. These formulas give \eqref{formulone}. 

Both sums in the right-hand side of \eqref{formulone} are determined because, roughly speaking, the first sum corresponds to the absolutely continuous part of the derivative and the second sum to its singular part (see \cite{Bi06} for the details).
\end{proof}

\begin{proof}[Proof of Theorem~\ref{poli}] 
Let $F_i$, $\nu_i$ and $I_p$ be as in the statement of Lemma~\ref{utile}.  Consider the distribution defined by the second sum in \eqref{formulone}. This distribution determines its support, which we denote by $S(A,w)$,  and determines
\begin{equation}\label{parte_singolare_dist}
d(x):=\sum_{(i,j)\in I_{p}} w\cdot\nu_i\ w\cdot\nu_j\ 
\la_{1}(F_i\cap(F_j+x)), 
\end{equation}
for $\la_{1}$-a.e. $x\in S(A,w)$. Note that $S(A,w)\subset\cup_{(i,j)\in I_{p}:\nu_i\cdot w\neq0}(F_i-F_j)$.
Choose any $i\in\{1,\dots,m\}$ and let $I_{\nu_i}=\{j\in1,\dots,m:\nu_j=\pm\nu_i\}$. We recall that $\rotpih\nu_i$ is a rotation of $\nu_i$ by $\pi/2$. Then, for any $h>0$  sufficiently small, we have
\[
d(h\,\rotpih \nu_i)=
\sum_{j\in I_{\nu_i}} (w\cdot\nu_i)^2\la_{1}(F_j\cap(F_j+h\ \rotpih \nu_i))=
(w\cdot\nu_i)^2
\sum_{j\in I_{\nu_i}}(\la_1(F_j)-h).
\]
Choose $w$ so that $w\cdot\nu_i\neq0$. Since the previous function is different from $0$, $S(A,w)$ contains a segment containing $o$ and parallel to $F_i$. Moreover, we have 
\[
\frac{\pa d}{\pa\rotpih\nu_i}(0)=-(w\cdot\nu_i)^2|I_{\nu_i}|,
\]
and this formula provides the number of edges of $A$ parallel to $F_i$.

The set $K$ is a convex polygon, because $D\,K$ coincides with $\supp g_A$, which is a polygon. Since $g_A=g_K$ the distribution considered above has the same features as the corresponding one for a convex polygon. This implies the following consequences.
\begin{enumerate}[{C}1]
\item \label{c1} The  number of edges of $A$ parallel to $F_i$ is at most two.
\item \label{c2} We have 
\begin{equation*}
S(A,w)\subset(\pa\, \supp g_A)\cup\left(\cup_{i:\nu_i\cdot w\neq0} \nu_i^\perp\right),
\end{equation*}
and each segment in $S(A,w)$ is parallel to an edge of $\supp g_A$. 
\end{enumerate}
We need to prove only C\ref{c2}, since C\ref{c1} is obvious. Assume $A$ convex polygon. If $F_i$ and $F_j$ are parallel and $i\neq j$ then $\nu_i=-\nu_j$ and, by \eqref{facce_corpodiff}, $F_i-F_j$ is an edge  of $D A=\supp g_A$. Moreover, when $i=j$ $F_i-F_j$ is a segment contained in $\nu_i^\perp$. Since $D A$ has an edge orthogonal to $\nu_i$, for each $i$, by~\eqref{facce_corpodiff} with $u=\nu_i$, the property is proved.

%

To conclude the proof of the theorem we argue by contradiction and assume $A$ non-convex. We have $D A= D(\conv A)$, because otherwise $g_A\neq g_K$, by Proposition~\ref{insdiff}. Consider the edges of $A$ not contained in $\partial(\conv A)$. We may assume that they are $F_1,\dots, F_d$, for some $d<m$, 
Let us distinguish the following three cases.

\begin{itemize}

\item[1)] \emph{There exists an edge $F_k$,
$k\in \{1,\ldots,d\}$, which is not parallel to any edge
of $\conv A$.} In this case, if we choose $w$ so that $w\cdot\nu_k\neq0$,  $S(A,w)$ contains a segment parallel to $F_k$ which is not parallel to any edge of $D(\conv A)=D\,A=\supp g_A$. This contradicts C\ref{c2}. 

\item[2)] \emph{There exists an edge $F_k$,
$k\in \{1,\ldots,d\}$, parallel to exactly one edge
$M$ of $\conv A$.} Let us show that $M$ is an edge of 
$A$. We can write $M=[m_1, m_2]$, with $m_1,
m_2\in A$. Let $u$ be the unit outer normal vector to $\conv A$ at
$M$, i.e. $M=F(\conv A, u)$. The hypothesis defining this case implies 
that $F(\conv A, -u)$ is not an edge and is a single
point $m$.  Thus, $F(D(\conv A), u)=[m_1, m_2]-m$, by \eqref{facce_corpodiff}. As $DA=D(\conv A)$, \eqref{facce_corpodiff} implies $M=F(A,u)$.

Consider now a Cartesian coordinate system so that $(0,1)=u$. 
Clearly $M$ and $F_k$ are parallel to the $x$-axis. 
Among the edges $F_1,\ldots,F_d$ parallel to the $x$-axis consider those with the  smallest $y$-coordinate. Among these edges consider the edge $F_m$ with largest abscissa 
(see Fig.~\ref{figura_E}).
\begin{figure}[h]
\centering
\includegraphics[scale=0.5]{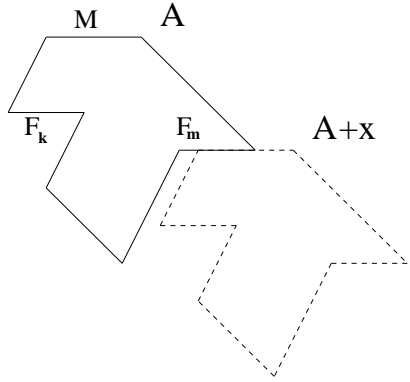}
\caption{The intersection of $A$ and $A+x$ for $x\in[x_1,x_0]$}\label{figura_E}
\end{figure}
Let $x_0$ be the translation which maps the point with smallest abscissa of $M$
to the point with largest abscissa of $F_m$, and let $x_1=x_0-(h,0)$ with $h>0$ sufficiently small. Note that $[x_1,x_0]\subset F_m-M$. We claim that, if $w\cdot(0,1)\neq0$, then   $[x_1,x_0]\subset S(A,w)$. 
Indeed, let $x\in[x_1,x_0]$ and consider the pairs of edges $F_i$ and $F_j$ of $A$ such that 
\begin{equation}\label{inter}
\la_1(F_i\cap(F_j+x))>0.
\end{equation}
If $F_i$ and $F_j$ are parallel to the $x$-axis, then we have necessarily $F_i=F_m$ and $F_j=M$, by the choice of $x_0$. If $F_i$ and $F_j$ are not parallel to the $x$-axis, then \eqref{inter} is false except possibly for finitely many $x\in[x_1,x_0]$. Therefore,  $\la_1$-a.e.\ in $[x_1,x_0]$ we have $d(x)=\pm w\cdot (0,1)\  \la_1(F_m\cap(M+x))\neq0$. This proves the claim.

The segment $[x_1,x_0]$ is not contained in a line through $o$, 
as $F_m$ is not aligned with $M$. Let us prove $[x_1,x_0]\not\subset\pa\supp g_A$.  Let $l_1$ and $l_2$ be the lines parallel to the $x$-axis supporting $A$, with $M\subset l_2$. The edges of $\supp g_A$ parallel to the $x$-axis are contained in $\pm(l_1-l_2)$. On the other hand, we have $[x_1,x_0]\subset F_m-M\not\subset l_1-l_2$, because $F_1\not\subset l_1$ ($l_1\cap A\subset l_1\cap\conv A=F(\conv A,(0,-1))$ and $F(\conv A,(0,-1))$ is not an edge, by the hypothesis defining this case). This proves $[x_1,x_0]\not\subset\pa\supp g_A$. 

These properties of $[x_1,x_0]$ contradict C\ref{c2}.

\item[3)] \emph{There exists an edge $F_k$,
$k\in \{1,\ldots,d\}$ parallel to a pair of antipodal parallel
edges $M$ and $N$ of $\conv A$.} Let us
show that at least one of the inequalities
\[
\lambda_1(M\cap\partial A)>0
\quad\text{and}\quad
\lambda_1(N\cap\partial A)>0
\]
holds. Let $u \in \mathbb{S}^1$ be such that
$M=F(\conv A, u)$ and $N=F(\conv A, -u)$ and assume that both inequalities are false. 
The geometric structure of $A$ implies that  both $F(A, u)$ and $F(A, -u)$
consist of a finite number of points. Consequently, $F(DA, u)$
consists of  a finite number of points, by \eqref{facce_corpodiff}, contradicting 
$DA= D(\conv A)$.

If exactly one of the previous inequalities holds, then the proof is concluded as in the previous case. If both  inequalities
hold, then $A$ has at least three edges orthogonal to $u$. This contradicts C\ref{c1}.

\end{itemize}
The above three cases complete all the possibilities. \end{proof}

\section{Non-convex sets with equal covariogram}\label{non-conv}

Gardner, Gronchi and Zong \cite{GGZ05}
presents  a pair of non-congruent non-convex polygons $P$ and $Q$ with equal covariogram.
\begin{figure}[h]
\centering
\includegraphics[scale=0.55]{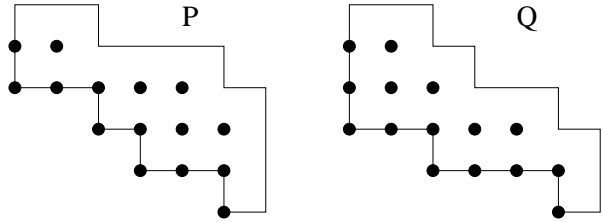}
\caption{Two non-congruent non-convex polygons with equal covariogram, which arise as animals of two homometric convex polyominoes (from \cite{GGZ05})}\label{figronchi}
\end{figure}
The polygons $P$ and $Q$ are the animals associated to two homometric convex
polyominoes consisting of fifteen points. 
We are interested in finding similar examples with minimal cardinality. Let us first prove that two animals have the same covariogram if and only if the corresponding polyominoes are homometric. The ``if'' part is proved, in a more general setting, in \cite{GGZ05}.

\begin{lemma}\label{poly_animals}
Let $A$ and $B$ be finite subsets of $\mathbb{Z}^n$ and let $\bar{A}=A+[0,1]^n$
and $\bar{B}=B+[0,1]^n$. Then $g_{\bar{A}}=g_{\bar{B}}$ if and only if $A$ and $B$ are homometric.
\end{lemma}
\begin{proof}
Let $Q=[0,1]^n$. 
\cite{GGZ05} proves the following formulas, valid for any $x\in\Real^n$,
\begin{equation}\label{cov_disc_cont}
g_{\bar{A}}(x)=\sum_{z\in \mathbb{Z}^n} |A\cap(A+z)| g_Q(z+x),\quad
g_{\bar{B}}(x)=\sum_{z\in \mathbb{Z}^n} |B\cap(B+z)| g_Q(z+x).
\end{equation}
If $A$ and $B$ are homometric these formulas imply $g_{\bar A}=g_{\bar B}$. Assume now $g_{\bar A}=g_{\bar B}$ and choose $w\in \mathbb{Z}^n$. The support of $g_Q(\cdot-w)$ is $D Q+w=[-1,1]^n+w$. Since $\mathbb{Z}^n\cap\inte([-1,1]^n+w)=\{w\}$, we have $g_Q(z-w)=0$ for each $z\in\mathbb{Z}^n$, $z\neq w$. Thus $g_{\bar A}(-w)=g_{\bar B}(-w)$ and \eqref{cov_disc_cont} imply
\[
 |A\cap(A+w)| =|B\cap(B+w)|.
\]
Since $|A\cap(A+w)| =|B\cap(B+w)|=0$ when $w\notin\mathbb{Z}^n$, the previous identity implies $A$ and $B$ homometric.
\end{proof}

The following proposition is known in the literature on homometric sets (see Rosenblatt and Seymour
\cite{RS}). It provides a method to construct pairs of homometric sets in any dimension. In some cases the obtained sets are polyominoes.

\begin{proposition}\label{hom}
Let $A$ and $B$ be subsets of $\mathbb{Z}^n$. Assume that each
point of $A+B$ (and of $A-B$) can be written in an unique way as sum
of a point of $A$ and of a point of $B$ (of -$B$, respectively).
Then $A+B$ and $A-B$ are homometric sets.
\end{proposition}

\begin{figure}[h]
\centering
\includegraphics[scale=0.55]{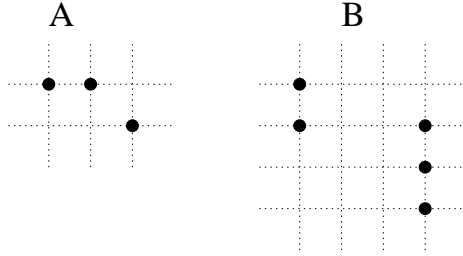}
\caption{The polyominoes in Fig.~\ref{figronchi} are equal to $A+B$ and $A-B$.}\label{figronchini}
\end{figure}

\begin{figure}[h]
\centering
\includegraphics[scale=0.55]{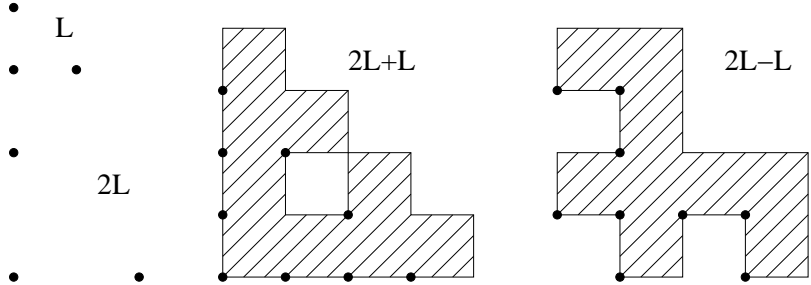}
\caption{Two non-congruent homometric polyominoes made of $9$ points, and the associated animals.}\label{figura_D}
\end{figure}

The example provided in \cite{GGZ05} can be obtained using this
construction. The pair of homometric polyominoes in Fig.
\ref{figronchi} can be written as $A+B$ and $A-B$, where $A$
and $B$ are the finite sets in Fig.~\ref{figronchini}.
Consider now the two sets of three points, $L$ and $2L$, in Fig.
\ref{figura_D}, and the two sets $2L+L$ and $2L-L$. These two sets
are homometric  polyominoes made of nine points. The corresponding
animals are non-congruent.

Another pair of animals made of nine squares which are not translations or reflections (with respect to a point) of each other is presented in \cite[Fig.~1]{DGN05}. The corresponding polyominoes are convex and one animal is the reflection of the other with respect to a line.

%

\begin{proof}[Proof of Theorem~\ref{polyominoes}] Let us consider two polyominoes $A,
B\subset\mathbb{Z}^2$ and the covariograms of $\bar{A}=A+[0,1]^2$
and $\bar{B}=B+[0,1]^2$. Obviously, $g_{\bar{A}}=g_{\bar{B}}$ implies 
$D\,\bar{A}=D\,\bar{B}$ as
$D\,\bar{A}=\supp g_{\bar{A}}$ and
$D\,\bar{B}=\supp g_{\bar{B}}$. Thus, the widths of $\bar{A}$
and $\bar{B}$ in the coordinate directions are equal. This
implies that the minimum rectangle with edges parallel to the coordinate axes containing $\bar{A}$ has to be
equal to the one containing $\bar{B}$. 

Let us denote by
$dP(h\times b)$ the class of $d$-polyominoes  (polyominoes
consisting of $d$ points) $A$ such that the minimal rectangular
container of $A+[0,1]^2$ has height $h$ and basis $b$. 
Let us remark that commonly polyominoes are classified up to all the
symmetries with respect to the coordinates axes. Here, however, we
will classify polyominoes up to translations and reflections in a
point, i.e.\ we identify
two polyominoes in $dP(h\times b)$ if they are reflections or
translations of each other.

We consider now the $d$-polyominoes, for each $d=1,\ldots,8$. 
It suffices to consider in the proof only polyominoes in
$dP(h\times b)$, with $h\leqslant b$. Indeed, the polyominoes in
$dP(b\times h)$ are obtained from those in $dP(h\times b)$ by a
rotation of $\pi/2$, and, moreover, a polyomino in $dP(h\times
b)$ cannot have the same covariogram of one in $dP(b\times h)$,
unless $h=b$, for the reason explained above.

The case $d=1$ and $d=2$ are trivial because there exist only one
1-polyomino and only one 2-polyomino that belongs to $2P(1 \times
2)$.

The class $3P(1\times 3)$ contains one element, while $3P(2\times
2)$ contains two elements. The two polyominoes in $3P(2\times 2)$
cannot have the same covariogram as their difference bodies are
not equal.

For $d=4$ the only class $4P(h\times b)$ with more than one
element is $4P(2\times 3)$. The five sets in $4P(2\times 3)$ have
different difference bodies.

For $d=5$ there are six elements in $5P(2\times 4)$, three
elements in $5P(2\times 3)$ and twelve elements in $5P(3\times
3)$. None of these sets has  difference body equal to that of another
set in the same class.

The elements in $6P(h\times b)$, in $7P(h\times b)$ and in
$8P(h\times b)$ have been analysed using the simple algorithm described in the appendix. In the case of 6-polyominoes, 7-polyominoes and
8-polyominoes the algorithm stops without finding a pair of
homometric polyominoes. 
\end{proof}


\section{Appendix}

\begin{figure}[h]
\centering
\includegraphics[scale=0.36]{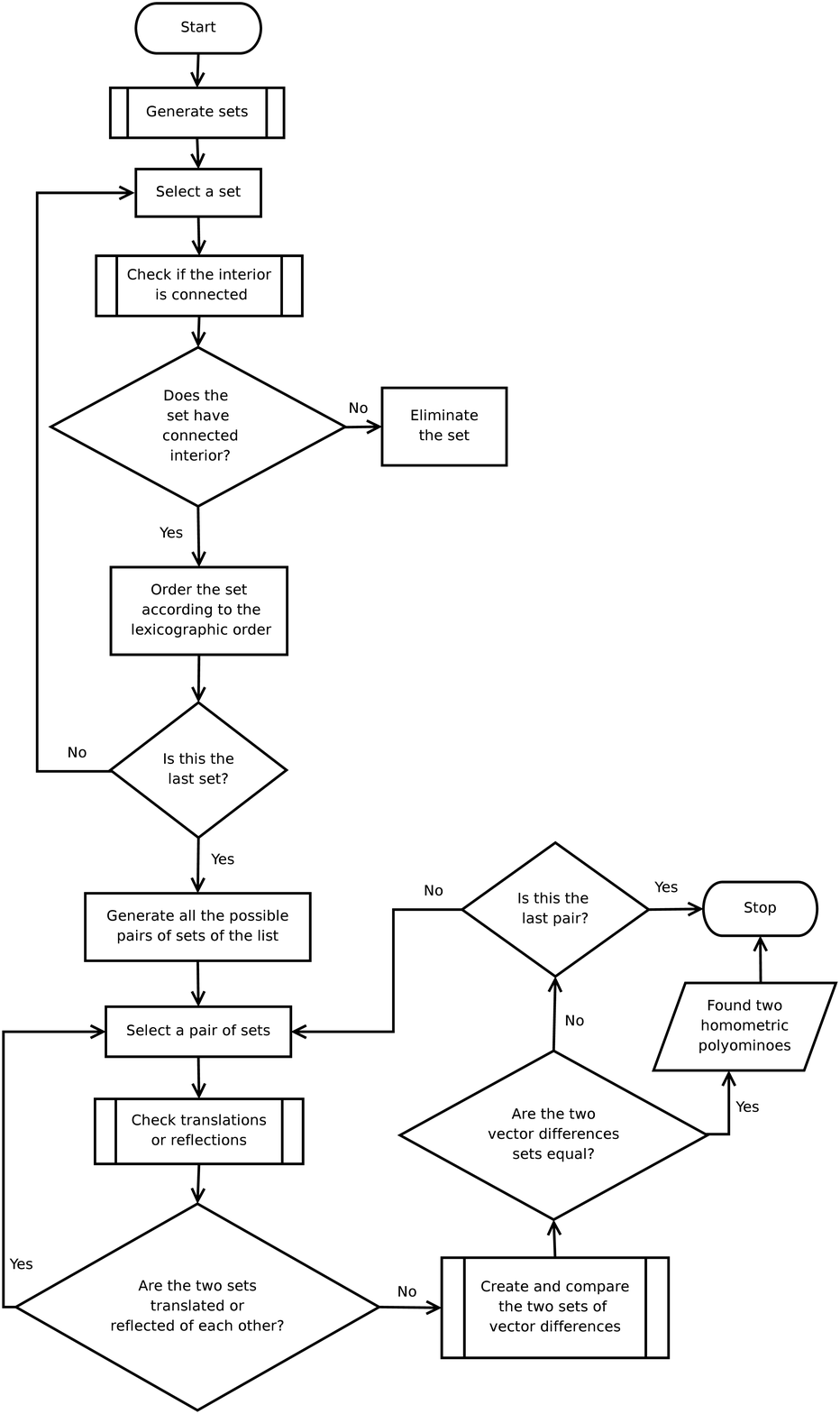}
\caption{The algorithm  used in the proof of Theorem~\ref{polyominoes}.}\label{diagramma}
\end{figure}

The diagram in Fig.~\ref{diagramma} describes the algorithm  used in the proof of Theorem~\ref{polyominoes}. We briefly explain what each subprogram does.

\begin{description}

\item[Generate sets] this function generates all
possible sets of eight (respectively seven, six) points of a grid
of $\mathbb{Z}^2$ with at most four (respectively four, three)
rows and eight (respectively seven, six) columns.
\item[Check if the interior is connected] this function chooses a
point $x_1$ of the selected set and constructs the  component containing the point. Successively it establishes if this  component
coincides with the whole set. It works with two lists
of points. At the beginning the first list $L_1$ contains only $x_1$, whereas  the second list $L_2$ contains all the other points of
the set. Among the points in $L_2$, the program transfers in $L_1$ those whose
distance from $x_1$ is unitary. Successively, the program considers the
second point in $L_1$ and repeats the process. The algorithm stops when it
has considered the last point in $L_1$. The set is connected if at the end
$L_2$ is empty.

\item[Check translations or reflections] this function
computes the vector differences of each point of the first set
$P_i$ with the corresponding (in the lexicographic order) point of
the other set, $P_j$. If all these differences are equal then the
two sets are translations of each other. If some of these differences are not equal,  then the function computes the
vector differences of each point of the first set with the
corresponding (in the lexicographic order) point of the second
set, previously reflected and ordered. If all these differences
are equal, the two sets are reflections of each other.
Otherwise $P_i$ and $P_j$ are not one translations or reflections of each other.
\item[Create and compare the two sets of vector
differences] this function generates for the pair $(P_i, P_j)$
the vector differences sets $DP_i$ and $DP_j$. Successively, it
orders $DP_i$ and $DP_j$ according to the lexicographic order and
compute the vector differences of each point of $DP_i$ with the
corresponding point of $DP_j$. If all these vectors are equal to
the null vector, then $P_i$ and $P_j$ are homometric. Otherwise
they are not homometric.
\end{description}

\end{document}